\documentclass[12pt, reqno]{amsart}
\usepackage{amssymb,amsmath}
\usepackage{amsmath, amsthm, amscd, amsfonts, amssymb, graphicx, color}
\usepackage{amsmath, amsthm, amssymb}
\usepackage{amsfonts}
\usepackage{graphicx}

\newtheorem{theorem}{Theorem}[section]
\newtheorem{lemma}{Lemma}[section]
\newtheorem{corollary}{Corollary}[section]
\newtheorem{definition}{Definition}[section]

\newtheorem{proposition}{Proposition}[section]

\newcommand{\M}{\mathbb{M}_{n}}
\newcommand{\MMM}{\mathbb{M}^s_{n, \alpha}}

\begin{document}

\title[General Numerical Radius for Products of Sectorial Matrices]{General Numerical Radius for Products of Sectorial Matrices}

\author[M. Alakhrass]{Mohammad Alakhrass}	

\address{ Department of Mathematics, University of Sharjah , Sharjah 27272, UAE.}
\email{\textcolor[rgb]{0.00,0.00,0.84}{malakhrass@sharjah.ac.ae}}



\begin{abstract}
In this paper, we investigate the generalized numerical radius $\omega_N$, associated with a matrix norm $N$ defined by \( \omega_N(X) = \sup_{\theta \in \mathbb{R}} N(\operatorname{Re}(e^{i\theta}X)) \).    
We focus on matrices whose numerical ranges are contained in sectors of the complex plane (sectorial matrices) and derive upper bounds for \( \omega_N(XY) \) and \( \omega_N(X \circ Y) \) for such matrices \( X \) and \( Y \). 
Our results generalize and refine well known numerical radius inequalities. Several known inequalities for \( \omega(X) \) are recovered as special cases.
\end{abstract}

\maketitle

\section{Introduction}

Let $\M$ denote the algebra of all $n \times n$ complex matrices. The Hadamard product of two matrices $X=[x_{i,j}]$ and $Y=[y_{i,j}]$ in $\M$ is defined as the matrix 
$$
X \circ Y = [x_{i,j}y_{i,j}].
$$
A matrix $X \in \M$ can be expressed as
\begin{equation}\label{cartesian decomposition}
X=\operatorname{Re}(X)+i\operatorname{Im}(X),
\end{equation}
where $\operatorname{Re}(X)$ and $\operatorname{Im}(X)$ are Hermitian matrices defined as 
$\operatorname{Re}(X) = \frac{X + X^*}{2}$ and $\operatorname{Im}(X) = \frac{X - X^*}{2i}$, respectively. This decomposition is commonly referred to as the Cartesian decomposition of $X$.

A matrix $X$ is called accretive (respectively, dissipative) if, in its cartesian decomposition \eqref{cartesian decomposition},  the matrix $\operatorname{Re}(X)$ (respectively, $\operatorname{Im}(X)$) is positive definite. If both $\operatorname{Re}(X)$ and $\operatorname{Im}(X)$ are positive definite, the matrix $X$ is referred to as accretive-dissipative.

A norm $\|\cdot\|$ on $\mathbb{M}_n$ is said to be unitarily invariant if it satisfies the property $\|UXV\|$ for all $X\in\mathbb{M}_n$ and all unitary matrices $U,V\in\mathbb{M}_n.$ An example of such a norm is the operator norm, defined as 
$$
\|X\|=\sup_{\|x\|=1}\|Xx\|.
$$

For a matrix $X \in \M$, we define its numerical range as follows:  
$$
W(X)=\{ \langle {Xx,x} \rangle: x \in \mathbb{C}^n, ||x||=1 \},
$$
where $\langle \cdot, \cdot \rangle$ denotes the standard inner product on $\mathbb{C}^n$, and $\|\cdot\|$ represents the Euclidean norm on $\mathbb{C}^n$. 
It is well known that $W(X)$ is a compact and convex subset of the complex plane.

Let $\alpha \in [0, \pi/2)$ and let $S_{\alpha}$ denote the sector of complex numbers defined as  
$$
S_{\alpha}=\{z=x+i y \in \mathbb{C}: x > 0, |y| \leq \tan(\alpha) x  \}.
$$
If the numerical range of a matrix $X$ is subset of $S_{\alpha}$, then $X$ is called sectorial. 
The smallest $\alpha$ for which this condition holds is referred to as the index of sectoriality. Let $\MMM, \alpha \in [0, \pi/2) $ represent the class of all $n \times n$ matrices $X$ such that $W(zX) \subset S_{\alpha}$, where $z$ is a complex number  
satisfying $|z|=1$. Let 
$$
\omega(X)=\sup  \{ |z|: z \in W(X)\}.
$$
The function $\omega(X)$ is usually refer to as the numerical radius of $X$. 

It is well-established that $\omega(\cdot)$ constitutes a norm on $\M$,   which satisfies the following inequalities:                                                                                     
 
\begin{equation}\label{A}
\frac{1}{2} \|X\|  \leq \omega(X) \leq   \| X \|; \quad  \forall X \in \M.
\end{equation}
Furthermore, if $X \in \mathbb{M}_n$ is a normal matrix, then $\omega(X) = \|X\|$, demonstrating that the bounds in \eqref{A} are sharp.

Clearly, $\omega( \cdot )$ represents a weakly unitarily invariant norm on $\M$; meaning it satisfies the condition
$\omega(UXU) = \omega(X)$ for all $X\in\mathbb{M}_n$   and any unitary matrix $U \in \mathbb{M}_n.$ However, $\omega( \cdot )$ is not necessarily a unitarily invariant norm.

It was established in \cite{yama} that for every $X \in \M$, the following holds:
$$
\omega(X)=\sup_{\theta\in \mathbb{R}}  ||Re (e^{i\theta}X)||. 
$$

Motivated by this result, the authors of \cite{ak} defined the generalized numerical radius as follows:

\begin{definition}\label{ND}
Let $ X \in \M$ and let $N$ be a norm on $\M$. The generalized numerical radius of $X$, induced by $N$, is defined as follows: 
$$
\omega_N(X)= \sup_{\theta\in \mathbb{R}} N( Re(e^{i\theta}X)).
$$
\end{definition}

They established several properties similar to those of $\omega(\cdot)$. 
Notably, they demonstrated that $\omega_N(\cdot)$ defines a norm on $\M$. 
Moreover, they showed that if $N(\cdot)$ is weakly unitarily invariant, then $\omega_N(\cdot)$ retains this invariance.

We recall that a norm $N$ on $\M$ is called multiplicative if satisfies:   
$$
N(X  Y) \leq N(X)  N(Y), \quad \text{for all $X ,Y \in \M$}, 
$$
and it is called Hadamard multiplicative if it satisfies: 
$$
N(X  \circ  Y) \leq N(X)   N(Y), \quad \text{for all $X ,Y \in \M$}. 
$$

The norm $\omega(\cdot)$ is neither multiplicative nor Hadamard multiplicative, in general. However, it is relatively straightforward to show that 

\begin{equation}\label{B}
\omega(X  Y) \leq  4 \ \omega(X)  \omega(Y); \quad  \forall X, Y \in \M,
\end{equation}
where the constant $4$ is optimal in \eqref{B}. Additionally, we have, 

\begin{equation}\label{C}
\omega(X   \circ   Y) \leq 2 \ \omega(X) \  \omega(Y),
\end{equation}
where the constant $2$ is optimal in \eqref{C}. See \cite[p. 73]{Horn and Jonson-Topics-book-1991}.

By selecting specific forms of $X$ and $Y$, it is possible to derive sharper estimates than those presented in \eqref{B} and \eqref{C}. 
For example, if $A, X \in \M$ and $A = [a_{ij}]$, then the following inequality holds:
\begin{equation}\label{I}
\omega(A \circ X) \leq   \left( \max_j a_{jj} \right) \, \, \omega(X).
\end{equation}
For further details, refer to \cite[Corollary 4]{Ando-1991} and \cite[Proposition 4.1]{Gau-Wu-2016}.

If $X=[x_{ij}] \in \mathbb{M}^s_{n, \alpha_1}$ and $Y=[y_{ij}] \in \mathbb{M}^s_{n, \alpha_2}$. Then, the following inequalities hold:

\begin{equation}\label{II}
\omega(XY) \leq \sec(\alpha_1)\sec(\alpha_2) \omega(X)\omega(Y),
\end{equation}

\begin{equation}\label{III}
\omega(X \circ Y) \leq \sec(\alpha_1)\sec(\alpha_2)  \omega(X)\omega(Y),
\end{equation}
and 
\begin{equation} \label{VI}
\omega(X \circ Y) \leq \sec(\alpha_1)\sec(\alpha_2)
\min \left\{\max_{j} | x_{jj} |  \, \,  \omega(Y), \, \max_{j} | y_{jj} |  \, \,  \omega(X)    \right \}.
\end{equation}
For further details, see \cite{Alakhrass-2023}.

In the sequel, unless stated otherwise,  let $N$ be a multiplicative, unitarily invariant, and self-adjoint norm, and let $\omega_N$ 
be defined as in the Definition \ref{ND}. 

This note aims to investigate the general numerical radius $\omega_N$ when applied to the product and the Hadamard product of sectorial matrices. 
Several properties will be established, including Propositions \ref{P1}, \ref{P2}, and \ref{P3}. 
Furthermore, the inequalities \eqref{I}, \eqref{II}, \eqref{III}, and \eqref{VI} are extended to the generalized numerical radius \( \omega_N \).

\section{Main Results}
We begin this section by presenting the following three lemmas. The proofs of these lemmas are available in 
\cite{Alakhrass-2021},  \cite{Alakhrass-Lieb} and \cite{Alakhrass-A note in sectorial matrices}. 

\begin{lemma}\label{L1}
Let $X \in \M$. If the numerical range of $X$ is a subset of $S_{\alpha}$, then   
$$
N \left( X \right) \leq  \sec(\alpha) N \left( Re (X)\right) .
$$
\end{lemma}

\begin{lemma}\label{L2}
Let $X \in \M$. If the numerical range of $X$ is a subset of $S_{\alpha}$, then 
$$
\left(
\begin{array}{cc}
\tan(\alpha) Re(X) & Im(X) \\
 Im(X) & \tan(\alpha) Re(X) \\
  \end{array}
\right) \geq 0. 
$$
\end{lemma}

\begin{lemma}\label{L3}
Let $X \in \M$. If the numerical range of $X$ is a subset of $S_{\alpha}$, then 

$$
\left(
  \begin{array}{cc}
   \sec(\alpha) Re(X) & X \\
    X^* & \sec(\alpha) Re(X) \\
  \end{array}
\right)\geq 0.
$$
\end{lemma}

In the following propositions, we present some basic properties of $\omega_N(\cdot)$ that will be used in the proofs.

\begin{proposition}\label{P1}
Let $X \in \M$. Then 
$$
\omega_N \left(Re(X) \right) \leq \omega_N \left( X \right).
$$ 
\end{proposition}

\begin{proof}
First of all, we remark that the fact that $Re(X^*)=Re(X)$ implies 
$\omega_N(X^*)= \omega_N(X)$. 
Now, observe that 
\begin{align}
\omega_N \left(Re(X) \right) &= \omega_N \left(\frac{X + X^*}{2} \right)  \notag \\
&= \frac{1}{2} \left( \omega_N \left( X+ X^* \right) \right) \notag \\
&\leq  \frac{1}{2} \left( \omega_N \left( X \right)  + \omega_N \left(X^* \right)  \right) \notag \\
&=  \frac{1}{2} \left( \omega_N \left( X \right)  + \omega_N \left(X \right)  \right) \notag \\
&=  \omega_N \left( X \right).  \notag
\end{align}
\end{proof}

\begin{proposition}\label{P2}
Let $X \in \M$. If the numerical range of $X$ is a subset of $S_{\alpha}$, then  
$$
\omega_N(Im(X)) \leq \tan(\alpha)  \omega_N(Re(X)).  
$$ 
\end{proposition}
\begin{proof}
Let $X=A + i B$ be the cartesian decomposition of $X$.
Since $W(X) \subset S_{\alpha}$, by Lemma \ref{L2}, we have   
$$
\left(
\begin{array}{cc}
\tan(\alpha) A & B \\
 B & \tan(\alpha) A \\
  \end{array}
\right) \geq 0.
$$
Thus,  
\begin{equation}\label{NN}
N(B) \leq \tan(\alpha) N(A).
\end{equation} 
Since $B$ is Hermitian, we have $\operatorname{Re}( e^{i \theta} B )=\cos(\theta) B.$ Therefore,  
\begin{align}
\omega_N(B) &= \sup_{\theta \in \mathbb{R} } N( \operatorname{Re}( e^{i \theta} B ) ) \notag  \\ 
& = \sup_{\theta \in \mathbb{R}} |\cos(\theta)| N(  B ) \notag  \\
&=N(B)  \notag  \\ 
& \leq \tan(\alpha)  N(A) \, (\text{by \eqref{NN}})\notag  \\   
& = \tan(\alpha)  \omega_N(A) \, (\text{since $A$ is Hermitian}) \notag  \\   
& = \tan(\alpha)  \omega_N(\operatorname{Re}(X)).\notag    
\end{align}
\end{proof}

\begin{proposition}\label{P3}
Let $X \in \M$. If the numerical range of $X$ is a subset of $S_{\alpha}$, then  
$$
\omega_N(X) \leq \sec(\alpha) \  \omega_N(\operatorname{Re}(X)).  
$$  
\end{proposition}
\begin{proof}
First, observe that by Lemma \ref{L1}, we have 
\begin{equation}\label{NN2}
N(X) \leq \sec(\alpha) N(\operatorname{Re}(X)).
\end{equation}
Now,
\begin{align}
\omega_N(X) &=\sup_{\theta \in \mathbb{R} } N( Re( e^{i \theta} X ) ) \notag  \\ 
&=\sup_{\theta \in \mathbb{R} } \frac{1}{2} N( e^{i \theta} X + e^{-i \theta} X^*  ) \notag  \\ 
&\leq \sup_{\theta \in \mathbb{R} } \frac{1}{2} \left( N( e^{i \theta} X) + N(e^{-i \theta} X^*)  \right) \notag  \\
&= \sup_{\theta \in \mathbb{R} } \frac{1}{2} \left( N( e^{i \theta} X) + N( e^{i \theta} X)  \right)= N( X)  \label{SA} \\ 
& \leq \sec(\alpha) N(\operatorname{Re}(X))  \, (\text{by \eqref{NN2} }) \notag  \\   
&= \sec(\alpha) \omega_N(\operatorname{Re}(X)). \notag  
\end{align}
\end{proof}

\subsection{Inequalities for Products of Sectorial Matrices}

The main result of this subsection is as follows:
\begin{theorem} \label{Th1}
Let $X\in \mathbb{M}^s_{n, \alpha_1}$ and $Y \in \mathbb{M}^s_{n, \alpha_2}$, with $\alpha_2, \alpha_2 \in [0, \pi/2)$. Then
$$
\omega_N(XY) \leq \sec(\alpha_1)\sec(\alpha_2) \omega_N(X)\omega_N(Y).
$$
\end{theorem}\label{Th-1}
\begin{proof}
Since $X \in \mathbb{M}^s_{n, \alpha_1}$ and $Y \in \mathbb{M}^s_{n, \alpha_2}$, we can find two complex numbers $z$ and 
$w $ such that $|z|=|w|=1$ where the numerical ranges  of $zX$ and $wY$ are subsets of $S_{\alpha_1}, S_{\alpha_2}$, respectively. Furthermore, we have 
\begin{align}
\omega_N(XY)  & \leq N \left(XY \right) \,  \  ( \text{by \eqref{SA}})  \notag \\
& \leq N \left(X \right)  \, \, N \left(Y \right)  \notag  \\
& =N \left( z X \right)  \, \, N \left(  w Y\right) \notag  \\
& \leq \sec(\alpha_1)\sec(\alpha_2)  N \left( Re \left(z X \right)  \right) N \left( Re \left( w  Y \right) \right)  \quad   \text{(by Lemma \ref{L1})}  \notag \\
& = \sec(\alpha_1)\sec(\alpha_2)  \omega_N( Re \left(z X \right)) \, \,   \omega_N \left(Re \left( w  Y \right) \right) \notag \\
& \text{(since both $Re \left(z X \right))$ and $ Re \left( w  Y \right)$ are Hermitian ) }  \notag \\
& \leq \sec(\alpha_1)\sec(\alpha_2) \omega_N \left(z X  \right)  \, \, \omega_N(w  Y)  \quad   \text{(by Proposition \ref{P1}) }  \notag \\
& = \sec(\alpha_1)\sec(\alpha_2) \omega_N \left(X  \right)  \, \, \omega_N(Y).   \notag
\end{align}
\end{proof}

As a direct result of Theorem \ref{Th1}, we derive the following corollaries.

\begin{corollary}
If $X,Y \in \MMM$, then
\begin{equation}\label{B2}
\omega_N(XY) \leq \sec^2(\alpha)\omega_N(X)\omega_N(Y).
\end{equation}

\begin{corollary}
If $X, Y \in \M$ are accretive-dissipative, then
$$
\omega_N(XY) \leq 2 \ \omega_N(X)    \omega_N(Y).
$$
\end{corollary}
\begin{proof}
Since $X, Y \in \M$ are accretive-dissipative, we have $X, Y \in \MMM$ with $\alpha = \pi /4.$ Now, the result follows from the above theorem,  \ref{Th1}. 
\end{proof}

\end{corollary}

\begin{corollary}
Let $X_j\in \mathbb{M}^s_{n, \alpha_j}, j=1,2,3,...,m$. Then
$$
\omega_N\left( X_1 X_2 ... X_m \right) \leq \left( \prod_{j=1}^{m} \sec(\alpha_j) \right)  \omega_N(X_1)\omega_N(X_2)... \omega_N(X_m) .
$$
\end{corollary}

\begin{corollary}
If $X_1, X_2, ..., X_m \in \MMM$, then
$$
\omega_N \left( X_1 X_2 ... X_m \right) \leq  \sec^m(\alpha) \omega_N(X_1)\omega_N(X_2)... \omega_N(X_m).
$$
\end{corollary}

\begin{corollary}
Suppose that $X_1, X_2, ..., X_m \in \M$ such that each one of them is accretive and dissipative, then
$$
\omega_N \left( X_1  X_1 ...  X_m \right) \leq 2^{m/2} \, \omega_N(X_1)\omega_N(X_m)... \omega_N(X_m).
$$
\end{corollary}

\subsection{Inequalities for Hadamard Products of Sectorial Matrices}
We begin this subsection with the following auxiliary lemma:

\begin{lemma}\label{L5}
\cite[Theorem 5.5.7]{Horn and Jonson-Topics-book-1991}. 
If $N$ is a unitary invariant norm on $\M$ , then it is Hadamard sub multiplicative if and only if it sub multiplicative. 
\end{lemma}
The generalized numerical radius $\omega_N$ is not, in general, Hadamard multiplicative. However, in a special case, we can prove the following result.

\begin{theorem}\label{H2}
Let $X,Y \in \M$ such that such that at least one of $X$ or $Y$ is Hermitian. Then 
$$
\omega_N(X \circ Y) \leq \omega_N(X)\omega_N(Y).
$$
\end{theorem}
\begin{proof}
Suppose $X,Y \in \M$ such that $Y$ is Hermitian. 
Then  
\begin{align}
\omega_N(X \circ Y)) &= \sup_{\theta} N(Re(e^{i \theta}(X \circ Y))) \\ \notag
&= \frac{1}{2} \sup_{\theta} N(    e^{i\theta}(X \circ Y) +e^{-i\theta}(X \circ Y)^*   ) \\ \notag
&= \frac{1}{2} \sup_{\theta} N(    e^{i\theta}X \circ Y +e^{-i\theta}X^* \circ Y   ) \\ \notag
&= \frac{1}{2} \sup_{\theta} N(  \left(  e^{i\theta}X +e^{-i\theta} X^* \right) \circ Y   ) \\ \notag
&=\sup_{\theta} N(  Re \left(  e^{i\theta}X \right) \circ Y   ) \\ \notag
& \leq \sup_{\theta} N(  Re \left(  e^{i\theta}X \right) ) N( Y   ) \\ \notag
& = \omega_N(X ) \omega_N(Y). \notag
\end{align}

\end{proof}

The main theorem of this subsection is stated as follows:

\begin{theorem}\label{H3}
Let $X\in \mathbb{M}^s_{n, \alpha_1}$ and $Y \in \mathbb{M}^s_{n, \alpha_2}$. Then, the following inequality holds:
\begin{equation}\label{E1}
\omega_N (X \circ Y) \leq \sec(\alpha_1)\sec(\alpha_2)  \omega_N(X)\omega_N(Y).
\end{equation}
\end{theorem}
\begin{proof}
Since $X\in \mathbb{M}^s_{n, \alpha_1}$ and $Y \in \mathbb{M}^s_{n, \alpha_2}$, there exist two complex numbers
$z, w \in \mathbb{C}$ such that $|z|=|y|=1$ where the numerical ranges of  $z X$ and $wY$ are subsets of $S_{\alpha_1}$ and $S_{\alpha_2}$, respectively. 
Therefore, by Lemma \ref{L3}, the following two blocks are positive semidefinite. 
$$
\left(
  \begin{array}{cc}
    \sec(\alpha_1)Re(zX) & zX \\
    \overline{z} X^* &  \sec(\alpha_1) Re(zX)\\
  \end{array}
\right),
\, \, \,
\left(
  \begin{array}{cc}
    \sec(\alpha_2) Re(wY) & wY \\
    \overline{w}Y^* &  \sec(\alpha_2) Re(wY) \\
  \end{array}
\right).
$$
Hence, there Hadamard products,   
$$
\left(
  \begin{array}{cc}
    \sec(\alpha_1)\sec(\alpha_2) Re(zX) \circ Re(wY) &zw \left( X \circ Y \right) \\
    \overline{zw} \left (X \circ Y \right)^* &  \sec(\alpha_1)\sec(\alpha_2) Re(zX) \circ Re(wY) \\
  \end{array}
\right), 
$$
is also positive semidefinite. This implies that 
\begin{align}
N \left(X \circ Y \right) &= N \left( zw \left( X \circ Y \right) \right) \notag\\
&\leq  \sec(\alpha_1)\sec(\alpha_2) N \left( Re(zX) \circ Re(wY) \right). \notag
\end{align}
Therefore, we have 
\begin{align}
\omega_N (X \circ Y ) & \leq N \left( X \circ Y \right) \notag \\
& \leq  \sec(\alpha_1)\sec(\alpha_2)  N\left( Re(zX) \circ Re(wY) \right)  \label{E2} \\
& \leq  \sec(\alpha_1)\sec(\alpha_2)  N\left( Re(zX) \right)  N \left(Re(wY) \right)   \quad (\text{by Lemma  \ref{L5}})\notag \\
& = \sec(\alpha_1)\sec(\alpha_2)   \omega_N( Re(zX) \omega_N(Re(wY))  \notag \\ 
& \text{(since both $Re(zX)$ and $Re(wY)$ are Hermitian)}  \notag \\
&\leq  \sec(\alpha_1)\sec(\alpha_2)  \omega_N( zX ) \omega_N( wY ) \quad (\text{by Proposition \ref{P1}})    \notag \\
&= \sec(\alpha_1)\sec(\alpha_2)  |z| \omega_N( X )  |w| \omega_N( Y ) \notag \\
&=  \sec(\alpha_1)\sec(\alpha_2)   \omega_N( X )  \omega_N( Y ). \notag
\end{align}
\end{proof}

As a consequence, we have the following corollary:

\begin{corollary}
If $X,Y \in \MMM$, then
\begin{equation}\label{C2}
\omega_N(X \circ  Y) \leq \sec^2(\alpha)\omega_N(X) \omega_N(Y).
\end{equation}
\end{corollary}

A straightforward modification of the arguments used in the proof of Theorem \ref{H3} yields the following Theorem.

\begin{theorem}
Let $X_j \in \mathbb{M}^s_{n, \alpha_j}, j=1,2,...,m$. Then
$$
\omega_N(X_1 \circ ... \circ X_m) \leq  \left( \prod_{j=1}^{m}  \sec(\alpha_j) \right)  \omega_N(X_1)\omega_N(X_2)...\omega_N(X_m). 
$$
\end{theorem}

Consequently, 

\begin{corollary}
If $X_1,..., X_m \in \M$ are accretive-dissipative, then

$$
\omega_N(X_1 \circ ... \circ X_m) \leq 2^{m/2}  \omega_N(X_1)\omega_N(X_2)...\omega_N(X_m).
$$
\end{corollary}

\subsection{Inequalities Involving Diagonal Entries}

We need the following lemma which can be found in \cite{Horn and Jonson-Topics-book-1991}.  

\begin{lemma} \label{L6} \cite[Theorem 5.5.19]{Horn and Jonson-Topics-book-1991}

Let $X,Y\in \mathbb{M}_n$ such that $Y=[y_{ij}] > 0$. Then 
$$
N(X \circ Y ) \leq \max_{i} y_{ii} \ N(X).
$$
\end{lemma}

Applying Lemma \ref{L6}, we can establish the following result.

\begin{lemma}\label{L7}
\label{hadamard}
Let $X,Y\in \mathbb{M}_n$ such that $Y > 0$. Then, 
$$
w_N( X \circ Y)\leq \max_i y_{ii} ~w_N(X). 
$$
\end{lemma}

\begin{proof}
Observe that, 
\begin{eqnarray*}
w_N(X \circ Y)&=& \sup_{\theta}N (Re[e^{i\theta}(X \circ Y)])\\
&=& \sup_{\theta} N(Re[(e^{i\theta}X)\circ Y])\\
&=& \sup_{\theta} N\left(\frac{(e^{i\theta}X)\circ Y+((e^{i\theta}X)\circ Y)^*}{2}\right)\\
&=& \sup_{\theta} N\left( \frac{(e^{i\theta}X)\circ Y+(e^{i\theta}X)^* \circ Y^*}{2} \right)\\
&=& \sup_{\theta} N\left( \frac{ \left(  ( e^{i\theta}X )+( e^{i\theta}X )^* \right) \circ Y}{2}\right)\\
&=& \sup_{\theta} N \left(Re[e^{i\theta}X] \circ Y \right)\\
&\leq& \sup_{\theta}N(Re[e^{i\theta}X]\max_{i}y_{ii}) \,  \  (\text{by Lemma \ref{L6}}) \\
&=& \max_{i}y_{ii} \ \sup_{\theta}N(Re[e^{i\theta}X])\\
&=& \max_i y_{ii} ~ w_N(X).
\end{eqnarray*}
\end{proof}

The result below estimates $\omega(X \circ Y)$ using the diagonal entries of the sectorial matrices $X$ and $Y$.
                                            
\begin{theorem} 
Suppose $X=[x_{ij}] \in \mathbb{M}^s_{n, \alpha_1}$ and $Y=[y_{ij}] \in \mathbb{M}^s_{n, \alpha_2}$. Then
$$
\omega_N(X \circ Y) \leq \sec(\alpha_1)\sec(\alpha_2)
\max_{j} | x_{jj} |  \, \,  \omega_N(Y),
$$
and 
$$
\omega_N(X \circ Y) \leq \sec(\alpha_1)\sec(\alpha_2)
\max_{j} | y_{jj} |  \, \,  \omega_N(X).
$$
\end{theorem}

\begin{proof}
By inequality \eqref{E2}, we have 
$$
\omega_N (X \circ Y ) \leq  \sec(\alpha_1) \sec(\alpha_2) \omega_N( Re(zX) \circ Re(wY)),
$$
where $z$ and $w$ are two complex numbers such that $|z|=|w|=1$.
Since $\operatorname{Re}(zX)$ is positive semidefinite, Lemma \ref{L7} implies that
$$
\omega_N( \operatorname{Re}(zX) \circ Re(wY)) \leq  \max_j \, \operatorname{Re}(zx_{jj}) \, \omega_N \left(  \operatorname{Re}(wY) \right).
$$
Now, we have

\begin{align}
\omega_N (X \circ Y ) & \leq \sec(\alpha_1)\sec(\alpha_2) \omega_N( \operatorname{Re}(zX) \circ \operatorname{Re}(wY))  \notag \\
&\leq  \sec(\alpha_1)\sec(\alpha_2)  \max_j \, \operatorname{Re}(zx_{jj}) \, \omega_N \left(  \operatorname{Re}(wY) \right) \notag \\
&\leq  \sec(\alpha_1)\sec(\alpha_2)  \max_j \, |z x_{jj}|  \, \omega_N \left(  wY \right) \notag \\
& =    \sec(\alpha_1)\sec(\alpha_2)  \max_j \, |x_{jj}| \, \omega_N \left(  Y \right). \notag
\end{align}

Hence,
$$
\omega_N (X \circ Y ) \leq \sec(\alpha_1)\sec(\alpha_2)   \max_j \, |x_{jj}| \, \omega_N \left(  Y \right).
$$
A similar argument implies the second inequality.
\end{proof}

\begin{corollary}
Let $X=[x_{ij}] \in \mathbb{M}^s_{n, \alpha_1}$ and $Y=[y_{ij}] \in \mathbb{M}^s_{n, \alpha_2}$. Then
$$
\omega_N(X \circ Y) \leq \sec(\alpha_1)\sec(\alpha_2)
\min \left \{ \max_{j} | x_{jj} |  \, \,  \omega_N(Y), \max_{j} | y_{jj}| \, \,  \omega_N(X) \right\}.
$$
\end{corollary}

\begin{corollary}
If $X=[x_{ij}], Y=[y_{ij}] \in \M$ are accretive-dissipative, then
$$
\omega_N(X \circ Y) \leq 2 \ \min \{ \max_{j} | x_{jj} |  \, \,  \omega_N(Y), \max_{j} | y_{jj}| \, \,  \omega(X) \}.
$$
\end{corollary}

An alternative upper bound for $\omega(X \circ Y)$, where $X$ and $Y$ are sectorial, can be derived as follows.

\begin{theorem}
Let $X \in \mathbb{M}^s_{n, \alpha_1}$ and $Y \in \mathbb{M}^s_{n, \alpha_2}$. Then
$$
\omega_N(X \circ Y)
\leq \min  \left \{(1+ \tan \alpha_1)  \omega_N (Re \, X ) \omega_N (Y), (1+ \tan \alpha_2)\omega_N (X)  \omega_N (Re Y ) \right\}.
$$
Consequently,
$$
\omega_N(X \circ Y) \leq (1+ \tan \alpha)\omega_N (X)  \omega_N (Y ),
$$
where $\alpha= \max \{ \alpha_1, \alpha_1\}$.
\end{theorem}

\begin{proof}
Since the second inequality follows from the first one, it is enough to prove the first inequality. Let $X=A + i B$ be the cartesian decomposition of $X$. 
Then
\begin{align}\label{zzz}
\omega_N(X \circ Y) & = \omega_N((A + i B) \circ Y)    \notag \\
& = \omega_N(A \circ Y  + i B \circ Y) \notag \\
& \leq \omega_N(A \circ Y) + \omega_N (B \circ Y) \notag \\
& \leq \omega_N(A) \omega (Y) + \omega_N (B) \omega_N (Y) \notag \quad \text{(by Theorem \ref{H2})}  \notag \\
& \leq \omega_N(A) \omega (Y) + \tan(\alpha_1) \,  \omega_N (A) \omega (Y) \notag \quad \text{(by Proposition \ref{P2})}  \notag \\
& =(1+ \tan \alpha_1)\omega_N (A) \omega_N (Y) \notag \\
& =    (1+ \tan \alpha_1)\omega_N (Re(X) ) \omega_N (Y). \notag
\end{align}

Hence,
\begin{equation}\label{F5}
\omega_N(X \circ Y) \leq (1+ \tan \alpha_1)  \omega_N (Re X ) \omega (Y).
\end{equation}

Similarly, one can show that
\begin{equation}\label{F6}
\omega_N(X \circ Y) \leq (1+ \tan \alpha_2)  \omega_N (Re Y ) \omega (X).
\end{equation}

The result follows by combining \eqref{F5} and \eqref{F6}.
\end{proof}

\begin{corollary}
Let $X, Y \in \MMM$. Then
$$
\omega_N(X \circ Y)
\leq (1+ \tan \alpha)   \min  \left \{(\omega_N (Re \, X ) \omega_N (Y), \omega_N (X)  \omega_N (Re Y ) \right\}.
$$
Consequently,
\begin{equation}\label{Z}
\omega_N(X \circ Y) \leq (1+ \tan \alpha)  \omega_N (X) \omega_N (Y).
\end{equation}
\end{corollary}

\textbf{Declarations}\\
\textbf{Ethics approval: } \\ 
Not applicable. 
\\
\textbf{Funding information:} \\
Not applicable. 
\\
\textbf{Data availability :} \\ 
No data sets were generated or analysed during the current study.
\\
\textbf{Competing interests:} \\
The authors declare no competing interests.

\end{document}